\documentclass[11pt]{amsart}

\usepackage{amsmath}
\usepackage{amssymb}
\usepackage{mathrsfs}
\usepackage[all]{xy}
\usepackage{comment}
\usepackage{color}
\usepackage[colorlinks,citecolor=blue,urlcolor=black,linkcolor=black]{hyperref}

\newtheorem{theorem}{Theorem}[section]

\newtheorem{corollary}[theorem]{Corollary}
\newtheorem{obs}[theorem]{Observation}

\theoremstyle{definition}
\newtheorem{definition}[theorem]{Definition}

\newtheorem{question}[theorem]{Question}

\theoremstyle{remark}

\newcount\skewfactor
\def\mathunderaccent#1#2 {\let\theaccent#1\skewfactor#2
\mathpalette\putaccentunder}
\def\putaccentunder#1#2{\oalign{$#1#2$\crcr\hidewidth
\vbox to.2ex{\hbox{$#1\skew\skewfactor\theaccent{}$}\vss}\hidewidth}}
\def\name{\mathunderaccent\tilde-3 }

\def\smallbox#1{\leavevmode\thinspace\hbox{\vrule\vtop{\vbox
   {\hrule\kern1pt\hbox{\vphantom{\tt/}\thinspace{\tt#1}\thinspace}}
   \kern1pt\hrule}\vrule}\thinspace}

\newcommand{\cf}{{\rm cf}}

\def\qedref#1{$\qed_{\reforiginal{#1}}$}

\title{Tiltan and Superclub}
\author{Shimon Garti}
\address{Einstein Institute of Mathematics,
 The Hebrew University of Jerusalem,
 Jerusalem 91904, Israel}
\email{shimon.garty@mail.huji.ac.il}

\author{Saharon Shelah}
\address{Institute of Mathematics
 The Hebrew University of Jerusalem
 Jerusalem 91904, Israel
 and  Department of Mathematics
 Rutgers University
 New Brunswick, NJ 08854, USA}
\email{shelah@math.huji.ac.il}
\urladdr{http://www.math.rutgers.edu/\char`\~shelah}

\thanks{Research supported by ISF grant no. 1838/19. This is publication 1179 in the second author's list}

\begin{document}
\let\labeloriginal\label
\let\reforiginal\ref
\def\ref#1{\reforiginal{#1}}
\def\label#1{\labeloriginal{#1}}

\begin{abstract}
We show that one can force superclub with an arbitrarily large value of ${\rm cov}(\mathscr{M})$.
We prove that the club principle is consistent with an arbitrarily large value of ${\rm add}(\mathscr{M})$.
We also prove that if $\kappa$ is regular then superclub at $\kappa^+$ implies $Q(\kappa^+,\kappa^+,\kappa^+)$. \\
Nous prouvons que superclub est consistant avec un valeur arbitrairement \'{e}lev\'{e}e de ${\rm cov}(\mathscr{M})$.
Nous prouvons que tr\`{e}fle est consistant avec un valeur arbitrairement \'{e}lev\'{e}e de ${\rm add}(\mathscr{M})$.
Nous prouvons aussi que superclub en $\kappa^+$ implique $Q(\kappa^+,\kappa^+,\kappa^+)$ si $\kappa$ est cardinal regulier.
\end{abstract}

\maketitle

\newpage

\section{Introduction}

The diamond principle $\Diamond_\kappa$ is a prediction principle discovered by Jensen in \cite{MR309729}.
Suppose that $\kappa$ is regular and uncountable.
A $\Diamond_\kappa$-sequence is a sequence of sets $(A_\alpha:\alpha\in\kappa)$ such that $A_\alpha\subseteq\alpha$ whenever $\alpha\in\kappa$ and for every $A\subseteq\kappa$ the set $S_A=\{\alpha\in\kappa:A\cap\alpha=A_\alpha\}$ is a stationary subset of $\kappa$.
We say that $\Diamond_\kappa$ holds iff there exists a $\Diamond_\kappa$-sequence.
Notice that $\Diamond_{\kappa^+}$ implies $2^\kappa=\kappa^+$.

The club principle $\clubsuit_\kappa$ is a weaker prediction principle, discovered by Ostaszewski in \cite{MR438292}.
Denote the set of limit ordinals of $\kappa$ by $\lim(\kappa)$.
We say that $(T_\delta:\delta\in\lim(\kappa))$ is a $\clubsuit_\kappa$-sequence iff every $T_\delta$ is an unbounded subset of $\delta$ and for every $A\in[\kappa]^\kappa$ the set $T_A=\{\delta\in\lim(\kappa):T_\delta\subseteq A\cap\delta\}$ is a stationary subset of $\kappa$.

The difference between $\Diamond_\kappa$ and $\clubsuit_\kappa$ is two-fold.
Firstly, $\Diamond_\kappa$ is based on equality while $\clubsuit_\kappa$ is based on inclusion.
Secondly, $\Diamond_\kappa$ predicts every element of $\mathcal{P}(\kappa)$ while $\clubsuit_\kappa$ predicts only elements of $[\kappa]^\kappa$.
The latter property is responsible for the fact that $\clubsuit_{\kappa^+}$ is consistent with $2^\kappa>\kappa^+$.
Since $\Diamond_{\kappa^+}$ implies $2^\kappa=\kappa^+$ one can see that $\clubsuit_{\kappa^+}$ is strictly weaker than $\Diamond_{\kappa^+}$.

Motivated by a stubborn open problem of Juhasz about the connection between $\clubsuit_{\aleph_1}$ and the existence of Suslin trees, Primavesi defined in \cite{primavesi} an intermediate principle, dubbed as \emph{superclub}.
We phrase the definition in the generalized form.

\begin{definition}
\label{defsuperclub} Superclub. \newline
Assume that $\kappa=\cf(\kappa)>\aleph_0$.
\begin{enumerate}
\item [$(\aleph)$] A superclub sequence at $\kappa$ is a sequence of sets $(S_\alpha:\alpha\in\lim(\kappa))$ such that $S_\alpha$ is an unbounded subset of $\alpha$ for each $\alpha$ and for every $A\in[\kappa]^\kappa$ there exists $B\in[A]^\kappa$ so that $S_B=\{\alpha\in\lim(\kappa):B\cap\alpha=S_\alpha\}$ is a stationary subset of $\kappa$.
\item [$(\beth)$] Superclub holds at $\kappa$ iff there exists a superclub sequence for $\kappa$.
\end{enumerate}
\end{definition}

Comparing superclub to diamond and club it is clear from the definition that diamond at $\kappa$ implies superclub at $\kappa$ which implies, in turn, the club principle at $\kappa$.
Both implications are irreversible.
It has been proved in \cite{MR3787522} that superclub at $\kappa^+$ implies Galvin's property at $\kappa^+$ while $\clubsuit_{\kappa^+}$ is consistent with the failure of Galvin's property at $\kappa^+$.
It has been proved in \cite{MR3687435} that superclub at $\aleph_1$ is consistent with $2^{\aleph_0}>\aleph_1$ while $\Diamond_{\aleph_1}$ implies $2^{\aleph_0}=\aleph_1$.

The forcing of \cite{MR3687435} is based on the existence of a strongly inaccessible cardinal in the ground model.
A natural question is whether superclub with a large continuum has some consistency strength.
We shall give a negative answer.
Specifically, in the model of \cite{MR1489304} the continuum hypothesis fails and superclub holds.
An immediate consequence is the consistency of superclub with a large value of some cardinal characteristics like the dominating number $\mathfrak{d}$.

A more interesting challenge is the bounding number $\mathfrak{b}$ and the category invariant ${\rm add}(\mathscr{M})$.
We do not know whether superclub is consistent with a large value of these characteristics, even if the desired value is just $\omega_2$.
However, we shall force club with an arbitrarily large value of $\mathfrak{b}$ and ${\rm cov}(\mathscr{M})$ simultaneously, thus we will have a large value of ${\rm add}(\mathscr{M})$ as well.

From the definition it seems that superclub has more affinities with the club principle.
The prediction of superclub applies only to $[\kappa]^\kappa$ and based on inclusion since if $B\cap\alpha=S_\alpha$ then $S_\alpha\subseteq A\cap\alpha$ and possibly $S_\alpha\neq A\cap\alpha$.
Hence in the above mentioned features which distinguish club from diamond, superclub behaves like club and unlike diamond.
But there is another feature in which superclub is similar to diamond.

If $A\subseteq\kappa$ and $S_A$ is the stationary set of guesses given by $\Diamond_\kappa$ then the corresponding elements of the diamond sequence are coherent in the following sense.
If $\gamma,\delta\in S_A$ and $\gamma<\delta$ then $A_\delta$ end-extends $A_\gamma$.
This is a consequence of the equality which forms the prediction of diamond sequences.
If one assumes only club at $\kappa$ then this coherence fades away.
Actually, one can choose a club sequence and then shrink every element in the sequence by taking subsets of order-type $\omega$, in which case $A_\delta$ never extends $A_\gamma$.

From this point of view, but only from this one, superclub is similar to diamond.
Given $A\in[\kappa]^\kappa$ and choosing $B\in[A]^\kappa$, the superclub sequence acts on $B$ as a diamond sequence with equality as the predicate for prediction.
Consequently, the elements of the superclub sequence are coherent, despite the fact that they are only included in the original set $A$.
In other words, a superclub sequence is a coherent club sequence.
This idea stands behind the proof of Galvin's property from superclub, and it will be exploited in the last section.

Our notation is (hopefully) standard.
We shall use the Jerusalem forcing notation, so $p\leq q$ means that $p$ is weaker than $q$.
Occasionally we shall use the word \emph{tiltan} instead of \emph{club} when refering to the club principle.
The paper contains three additional sections.
In the first part we force superclub with large continuum without an inaccessible cardinal.
Actually, ${\rm cov}(\mathscr{M})$ assumes the value of the continuum in the generic extension, hence $\mathfrak{d}=2^\omega$ as well.
In the second part we show that one can force club with ${\rm add}(\mathscr{M})=\mathfrak{c}$ where $\mathfrak{c}$ is arbitrarily large.
In the last section we deal with a graph-theoretic statement which follows from the continuum hypothesis and here we show that it follows merely from superclub.

We thank the anonymous referee of the paper for a careful reading of the article and for the helpful comments, all of them are integrated within the manuscript.
We also thank Thilo Weinert for a helpful conversation concerning Cicho\'n's diagram and the bounding number.
Finally, we thank J\"org Brendle and Lajos Soukup for several important comments and useful information.

\section{A takeaway theorem}

In this section we prove that the method of \cite{MR1489304} for enlarging the continuum while preserving club sequences from the ground model works equally well with respect to superclub sequences.
The main conclusion is that one can force superclub with an arbitrarily large value of the continuum without the assumption that there is an inaccessible cardinal in the ground model.
Another conclusion is that one can force superclub with a large value of ${\rm cov}(\mathscr{M})$, a result which will be useful in the next section.

Let us describe the basic forcing notion of \cite{MR1489304} in the spirit of Prikry forcing, namely we define the forcing order which is quite incomplete and we also define a pure order which is $\omega_1$-complete.

\begin{definition}
\label{deffss} The Fuchino-Shelah-Soukup forcing. \newline
Let $(\mathbb{P},\leq,\leq^*)$ be the following forcing notion.
\begin{enumerate}
\item [$(\aleph)$] A condition $p\in\mathbb{P}$ is a partial function from $\omega_2$ into $\{0,1\}$ such that $|{\rm dom}(p)|\leq\aleph_0$ and ${\rm dom}(p)\cap[\delta,\delta+\omega)$ is finite whenever $\delta$ is a limit ordinal of $\omega_2$.
\item [$(\beth)$] If $p,q\in\mathbb{P}$ then $u_{pq}=\{\delta\in\lim(\omega_2): \varnothing\neq{\rm dom}(p)\cap[\delta,\delta+\omega) ={\rm dom}(q)\cap[\delta,\delta+\omega)\}$.
\item [$(\gimel)$] If $p,q\in\mathbb{P}$ then $p\leq q$ iff $p\subseteq q$ and $u_{pq}$ is finite.
We define $p\leq^*q$ iff $p\subseteq q$ and $u_{pq}$ is empty.
\end{enumerate}
\end{definition}

The following statement is proved in \cite[Proposition 5.1.4]{primavesi}.
If $\name{f}$ is a $\mathbb{P}$-name, $p\in\mathbb{P}$ and $p$ forces that $\name{f}:\omega_1^V\rightarrow\omega_1^V$ is a function then one can find a set $A_{p,\name{f}}\in[\omega_1]^{\omega_1}\cap V$ and a function $g_{p,\name{f}}:A_{p,\name{f}}\rightarrow\omega_1$ such that $g_{p,\name{f}}\in V$ and for every ordinal $\eta\in\omega_1^V$ there is a condition $q_\eta\in\mathbb{P}$ such that $p\leq q_\eta$ and $q_\eta\Vdash g_{p,\name{f}}\upharpoonright(A_{p,\name{f}}\cap\eta)= \name{f}\upharpoonright(A_{p,\name{f}}\cap\eta)$.
Similar assertions appear in \cite{MR1489304}, see for example Lemma 3.9 there.

It follows from this statement that $\aleph_1$ is preserved, and moreover stationary subsets of $\aleph_1$ remain stationary in the generic extension as follows from \cite[Theorem 2.6]{MR1489304}.
To see that $\aleph_1$ is preserved consider $\name{h}:\omega_1^V\rightarrow\omega$ and a condition $p$ which forces this fact.
Let $g$ be $g_{p,\name{h}}$, so ${\rm rang}(g)\subseteq\omega$.
Since $A_{p,\name{h}}={\rm dom}(g)$ is uncountable, there are $\beta,\gamma\in A_{p,\name{h}}$ for which $g(\beta)=g(\gamma)$.
Now if $\delta\in\omega_1^V$ is sufficiently large then $q_\delta\geq p$ and $q_\delta\Vdash\name{h}(\beta)=g(\beta)=g(\gamma)=\name{h}(\gamma)$ so $\name{h}$ cannot be one to one and hence $\omega_1^V$ is not a countable ordinal in $V[G]$.

In the above definition the domain of any condition $p$ is contained in $\omega_2$, but one can replace $\omega_2$ by an arbitrarily large cardinal $\kappa$.
If $\kappa^{\aleph_0}=\kappa$ in the ground model then $2^\omega=\kappa$ in the generic extension.

\begin{theorem}
\label{thmaway} Let $V$ be a model of GCH and let $(S_\alpha:\alpha\in\lim(\omega_1))$ be a superclub sequence in $V$.
Let $\kappa\geq\omega_2$ be such that $\kappa^{\aleph_0}=\kappa$, let $\mathbb{P}$ be the Fuchino-Shelah-Soukup forcing based on $\kappa$, and let $G\subseteq\mathbb{P}$ be generic over $V$.
Then cardinals and cofinalities are preserved in $V[G]$, the continuum becomes $\kappa$ and $(S_\alpha:\alpha\in\lim(\omega_1))$ remains a superclub sequence in $V[G]$.
Hence superclub is consistent with an arbitrarily large value of ${\rm cov}(\mathscr{M})$.
\end{theorem}

\par\noindent\emph{Proof}. \newline
By the comment before the statement of the theorem we know that $\aleph_1$ is preserved.
Cardinals above $\aleph_1$ are preserved since $\mathbb{P}$ is $\aleph_2$-cc (here we use the fact that $V$ models \textsf{GCH}).
It is easy to infer from density arguments that $2^\omega=\kappa$ in $V[G]$, as the generic set $G$ adds a function from $\kappa$ into $\{0,1\}$ which can be sliced into $\kappa$ distinct functions from $\omega$ into $\{0,1\}$.

Suppose that $A\in[\omega_1]^{\omega_1}\cap V[G]$ and let $\name{f}$ be a name for the increasing enumeration of the elements of $A$.
Fix a condition $p$ which forces this fact.
Let $A_{p,\name{f}}$ and $g_{p,\name{f}}$ be as described ahead of the proof, so $A_{p,\name{f}}\in V$ and for every $\delta\in\omega_1^V$ one can find $q_\delta\geq p$ such that $q_\delta\Vdash{\rm rang}(g_{p,\name{f}}\upharpoonright A_{p,\name{f}})\subseteq A$.
Applying superclub in the ground model to $g_{p,\name{f}}''A_{p,\name{f}}$, let $B$ be an uncountable subset of $g_{p,\name{f}}''A_{p,\name{f}}$ on which the superclub sequence $(S_\alpha:\alpha\in\lim(\omega_1))$ acts like diamond.

Let $S_B=\{\alpha\in\omega_1:B\cap\alpha=S_\alpha\}$.
If $\alpha\in S_B$ then choose $\delta\in(\alpha,\omega_1)$ such that $S_\alpha\subseteq{\rm rang}(g_{p,\name{f}}\upharpoonright(A_{p,\name{f}}\cap\delta))$.
Let $q_\delta\geq p$ be so that $q_\delta\Vdash g_{p,\name{f}}\upharpoonright(A_{p,\name{f}}\cap\delta)= \name{f}\upharpoonright(A_{p,\name{f}}\cap\delta)$.
It follows that $q_\delta\Vdash B\cap\alpha=S_\alpha$.
Since this is true for every $\alpha\in S_B$ and the stationarity of $S_B$ is preserved, one can see that superclub holds in the generic extension.
It is easy to verify that ${\rm cov}(\mathscr{M})$ equals $\mathfrak{c}$ in the generic extension, see \cite{MR2279653} for details, so we are done.

\hfill \qedref{thmaway}

Remark that the main result of this section already shows that in order to force superclub with a large value of the continuum there is no need of an inaccessible cardinal in the ground model.
In the next section we will show that the above result can be strengthened if one replaces superclub by tiltan.
The idea is that the forcing notion of \cite{MR1489304} can be composed with Hechler forcing, thus making both ${\rm cov}(\mathscr{M})$ and $\mathfrak{b}$ large.
We do not know how to preserve superclub sequences from the ground model in this construction, but if one begins with a diamond sequence in the ground model then it remains a club sequence in the generic extension.

\newpage

\section{The additivity of the meager ideal}

In this section we prove the consistency of tiltan (the club principle) with an arbitrarily large value of ${\rm add}(\mathscr{M})$.
This confirms a conjecture of Brendle, see \cite[Conjecture 9.4]{MR2279653}.
The context of this result is nicely displayed by Cicho\'n's diagram:

\xymatrix{
& & {\rm cov}(\mathscr{N}) \ar[dd] & {\rm non}(\mathscr{M}) \ar[l] \ar[d] & {\rm cof}(\mathscr{M}) \ar[l] \ar[d] & {\rm cof}(\mathscr{N}) \ar[l] \ar[dd] \\
&&& \mathfrak{b} \ar[d] & \mathfrak{d} \ar[l] \ar[d] \\
& & {\rm add}(\mathscr{N}) & {\rm add}(\mathscr{M}) \ar[l] & {\rm cov}(\mathscr{M}) \ar[l] & {\rm non}(\mathscr{N}) \ar[l] \\
}

A result of Truss from \cite{MR727790} says that $\clubsuit_{\aleph_1}$ implies ${\rm add}(\mathscr{N})=\aleph_1$, so the smallest invariant in the above diagram becomes $\aleph_1$ under the club principle.
Nevertheless, we shall see that $\clubsuit_{\aleph_1}$ is consistent with an arbitrarily large value of ${\rm add}(\mathscr{M})$.
The idea is to increase both ${\rm cov}(\mathscr{M})$ and $\mathfrak{b}$ while keeping the club principle. By a result of Miller from \cite{MR613787} we know that ${\rm add}(\mathscr{M})=\min(\mathfrak{b},{\rm cov}(\mathscr{M}))$, so our theorem will follow.

A natural way to increase ${\rm cov}(\mathscr{M})$ is by using \cite{MR1489304}, as described in the previous section.
A natural way to increase $\mathfrak{b}$ is Hechler forcing from \cite{MR0360266}.
We need, therefore, to amalgamate these two forcing notions.
Moreover, we would like to do it while maintaining tiltan.
For this end, we shall use a strong version of $ccc$, called \emph{sweetness}.
The concept of sweet forcing notions comes from \cite{MR768264}, and a good background is contained in \cite{MR2076408}.
Let us recall the formal definition.

\begin{definition}
\label{defsweet} Sweet forcing notions. \newline
A forcing notion $\mathbb{P}$ is sweet iff there are a dense subset $\mathcal{D}$ of $\mathbb{P}$ and a sequence $\bar{E}=(E_n:n\in\omega)$ of equivalence relations on $\mathcal{D}$ such that the following requirements are met:
\begin{enumerate}
\item [$(a)$] Each $E_n$ is $\leq_{\mathbb{P}}$-directed and $\mathcal{D}/E_n$ is countable.
\item [$(b)$] $E_{n+1}\subseteq E_n$ for every $n\in\omega$.
\item [$(c)$] If $\{p_i:i\in\omega\}\subseteq\mathcal{D}$ and $r\in\mathcal{D}$ satisfies $p_iE_ir$ for every $i\in\omega$ then for every $n\in\omega$ there is a condition $q_n\geq r$ such that $q_nE_nr$ and $\forall i\geq n, p_i\leq q_n$.
\item [$(d)$] If $p,q\in\mathcal{D}$ and $n\in\omega$ then there is some $m\in\omega$ so that for every $r\in[p]_{E_m}$ one can find $t\in[q]_{E_n}$ such that $r\leq t$.
\end{enumerate}
\end{definition}

There is another concept of sweetness, defined by Stern in \cite{MR808816} and based on topological considerations.
For the formal definition of topological sweetness and the concept of an iterably sweet forcing notion we refer to \cite[Section 4]{MR2076408}.
Here we just mention the fact that if $\mathbb{P}$ is sweet and every pair of compatible conditions has a least upper bound then $\mathbb{P}$ is iterably sweet.
The forcing notion of \cite{MR1489304} mentioned in the previous section is sweet and also iterably sweet.
It has been proved in \cite[Section 7]{MR768264} that if $\mathbb{P}$ is sweet and $\name{\mathbb{D}}$ is a $\mathbb{P}$-name of Hechler forcing then $\mathbb{P}\ast\name{\mathbb{D}}$ is sweet as well.

In what follows we describe a class of forcing notions which generalize the forcing of \cite{MR1489304}.
Define $\mathcal{K}_0$ as the class of objects $K$ which are essentially an iteration of $ccc$ forcing notions.
So $\ell g(K)$ is an ordinal (the length of the iteration) and $\mathbb{P}_\alpha$ is a forcing notion for every $\alpha\leq\ell g(K)$.
If $\alpha\leq\beta\leq\ell g(K)$ then we require $\mathbb{P}_\alpha\lessdot\mathbb{P}_\beta$, and for every $\beta<\ell g(K)$ we have a $\mathbb{P}_\beta$-name $\name{\mathbb{Q}}_\beta$ of a $ccc$ forcing notion.

Suppose that $\alpha\leq\ell g(K)$.
A condition $p\in\mathbb{P}_\alpha$ is a function so that ${\rm dom}(p)\subseteq\alpha$ and $|{\rm dom}(p)|\leq\aleph_0$.
If $\beta\in{\rm dom}(p)$ then $p(\beta)$ is a $\mathbb{P}_\beta$-name of an element of $\name{\mathbb{Q}}_\beta$ and if $\beta<\alpha$ then $p\upharpoonright\beta\in\mathbb{P}_\beta$.
If $p,q\in\mathbb{P}_\alpha$ then $p\leq_{\mathbb{P}_\alpha}q$ iff ${\rm dom}(p)\subseteq{\rm dom}(q)$ and for every $\beta\in{\rm dom}(p)$ we have $p\upharpoonright\beta\Vdash_{\mathbb{P}_\beta} p(\beta)\leq_{\name{\mathbb{Q}}}q(\beta)$, and the set ${\rm diff}(p,q)=\{\beta\in{\rm dom}(p):p(\beta)\neq q(\beta)\}$ is finite.

Apart from the usual ordering $\leq_{\mathbb{P}_\alpha}$ we define the orderings $\leq^{\rm pr}_\alpha$ and $\leq^{\rm ap}_\alpha$ where pr stands for pure and ap stands for apure.
Formally, $\leq^{\rm pr}_\alpha=\{(p,q)\in\mathbb{P}_\alpha\times\mathbb{P}_\alpha: p\leq_{\mathbb{P}_\alpha}q\wedge q\upharpoonright{\rm dom}(p)=p\}$ and $\leq^{\rm ap}_\alpha=\{(p,q)\in\mathbb{P}_\alpha\times\mathbb{P}_\alpha: p\leq_{\mathbb{P}_\alpha}q\wedge {\rm dom}(p)={\rm dom}(q)\}$.
Notice that if $p\leq_{\mathbb{P}_\alpha}r$ then one can separate the extension by finding $q\in\mathbb{P}_\alpha$ so that $p\leq^{\rm pr}_\alpha q\leq^{\rm ap}_\alpha r$.
If $K\in\mathcal{K}_0$ and $\alpha=\ell g(K)$ then $\leq_K=\leq_{\mathbb{P}_\alpha}, \mathbb{P}_K$ denotes the whole iteration and so forth.

Assuming \textsf{CH} in the ground model, if $K\in\mathcal{K}_0$ then $\mathbb{P}_K$ is $\aleph_2$-cc.
The pure order is $\aleph_1$-complete, and moreover if $(p_n:n\in\omega)$ is $\leq_K^{\rm pr}$-increasing then $\bigcup_{n\in\omega}p_n$ is an upper bound, which will be called the \emph{canonical} upper bound.
One can verify that $\mathbb{P}_K$ is proper.
Assume that $(p_i:i\in\omega_1)$ is $\leq^{\rm pr}_K$-increasing and continuous, that is $p_j$ is the canonical upper bound of $(p_i:i\in j)$ whenever $j\in\omega_1$ is a limit ordinal.
Assume further that $p_{i+1}\leq^{\rm ap}_K q_{i+1}$ for every $i\in\omega_1$.
Then one can find $i<j<\omega_1$ such that $q_{i+1}\parallel q_{j+1}$.
It can be verified that the forcing notion of \cite{MR1489304} used in the previous section belongs to $\mathcal{K}_0$.
Let $\mathcal{K}$ be the subclass of $\mathcal{K}_0$ which consists of iterations $K$ which are Suslin $ccc$ (this means that the forcing and its order are $\Sigma^1_1$-definable) and iterably sweet.
We shall use this class in order to prove the following:

\begin{theorem}
\label{thmaddm} Assume that $\Diamond_{\aleph_1}$ holds in $V, K\in\mathcal{K}$ and $\name{\mathbb{Q}}_{K,\beta}$ is a name of Hechler forcing for every $\beta\leq\ell g(K)$.
Then $\clubsuit_{\aleph_1}$ holds in the generic extension by $\mathbb{P}_K$.
Concomitantly, $\mathfrak{b}=\ell g(K)$ in the generic extension by $\mathbb{P}_K$, provided that $\ell g(K)$ is a regular and uncountable cardinal.
\end{theorem}

\par\noindent\emph{Proof}. \newline
The fact that $\mathfrak{b}=\ell g(K)$ in the generic extension by $\mathbb{P}_K$ follows from the properties of Hechler forcing.
We must show, therefore, that $\clubsuit_{\aleph_1}$ holds in this generic extension.
Let $(A_\alpha:\alpha\in\omega_1)$ be a diamond sequence in $V$.
Since $\ell g(K)$ is typically large, we get the negation of \textsf{CH} in the generic extension and hence diamond fails.
However, we will show that $(A_\alpha:\alpha\in\omega_1)$ exemplifies $\clubsuit_{\aleph_1}$ in the generic extension.
Suppose that $\name{A}$ is a $\mathbb{P}_K$-name of an element of $[\omega_1]^{\aleph_1}$.
Fix a generic set $G\subseteq\mathbb{P}_K$ and a condition $p$ such that $p\Vdash\name{A}\in[\omega_1]^{\aleph_1}$.
By induction on $i\in\omega_1$ we choose a triple $(p_i,q_i,\alpha_i)$ such that the following requirements are met:
\begin{enumerate}
\item [$(a)$] $p_0=p$.
\item [$(b)$] $\alpha_i\in\omega_1$.
\item [$(c)$] $(p_j:j\leq i)$ is $\leq_K^{\rm pr}$-increasing and continuous.
\item [$(d)$] $p_i\leq_K q_i$.
\item [$(e)$] $q_i$ forces that the $i$th member of $\name{A}$ is $\check{\alpha}_i$.
\item [$(f)$] $p_{i+1}\leq_K^{\rm ap}q_i$.
\end{enumerate}
Let $\chi$ be a sufficiently large regular cardinal.
For every $\varepsilon\in\omega_1$ choose $N_\varepsilon\prec\mathcal{H}(\chi)$ such that each $N_\varepsilon$ is countable, $(N_\varepsilon:\varepsilon\in\omega_1)$ is increasing and continuous, $(N_\zeta:\zeta\leq\varepsilon)\in N_{\varepsilon+1}, K\in N_0$ and $((p_i,q_i,\alpha_i):i\in\omega_1)\in N_0$.
Fix an ordinal $\delta\in\omega_1$ so that $N_\delta\cap\omega_1=\delta$.

Let $w={\rm diff}(p_{\delta+1},q_\delta)=\{\alpha\in{\rm dom}(q_\delta): p_{\delta+1}(\alpha)\neq q_\delta(\alpha)\}$.
Choose a sequence $(\varepsilon_n:n\in\omega)$ such that $m<n\Rightarrow\varepsilon_m<\varepsilon_n, \delta=\bigcup_{n\in\omega}\varepsilon_n$ and the following property holds.
For every first order formula $\varphi(x,\bar{y})$ and every $\bar{a}\in{}^{\ell g(\bar{y})}N_\delta$ and every large enough $n\in\omega, \mathcal{H}(\chi)\models\varphi[q_\delta,\bar{a}]\equiv \varphi[q_{\varepsilon_n},\bar{a}]$.
For each $n\in\omega$ let $\delta_n=N_{\varepsilon_n}\cap\omega_1$, so $(\delta_n:n\in\omega)$ is increasing and $\delta=\bigcup_{n\in\omega}\delta_n$.
We define a condition $q$ as follows:
\begin{enumerate}
\item [$(\aleph)$] ${\rm dom}(q)=N_\delta\cap\ell g(K)$.
\item [$(\beth)$] If $\alpha\in{\rm dom}(q)\cap w$ then $q(\alpha)=p_\delta(\alpha)$.
\item [$(\gimel)$] If $\alpha\in{\rm dom}(q)-w$ then $q(\alpha)$ is a common upper bound of the conditions $(q_{\delta_n}:n\in\omega)$.
\end{enumerate}
The last item can be instantiated by sweetness.
Let $a=\{\alpha_{\delta_n}:n\in\omega\}$ and let $S$ be a stationary subset of $\omega_1$ for which $a=A_\alpha$ whenever $\alpha\in S$.
Here we use the fact that $(A_\alpha:\alpha\in\omega_1)$ is a diamond sequence in $V$.
Since $q_{\delta_n}\leq q$ for every $n\in\omega$ we see that $q\Vdash a=A_\alpha\subseteq\name{A}$ for every $\alpha\in S$, so we are done.

\hfill \qedref{thmaddm}

We can derive now the main result of this section:

\begin{corollary}
\label{coraddm} Tiltan and ${\rm add}(\mathscr{M})$. \newline
It is consistent that $\clubsuit_{\aleph_1}$ holds and ${\rm add}(\mathscr{M})$ is arbitrarily large.
\end{corollary}

\par\noindent\emph{Proof}. \newline
Assume $\Diamond_{\aleph_1}$ in $V$ and choose an arbitrarily large regular cardinal $\lambda$ which satisfies $\lambda^{\aleph_1}=\lambda$.
Let $K\in\mathcal{K}$ be such that $\lambda=\ell g(K)$ and $\name{\mathbb{Q}}_{K,\beta}$ is Hechler forcing for every $\beta\in\lambda$.
Let $G\subseteq\mathbb{P}_K$ be $V$-generic.

By the definition of $\mathcal{K}$ we use a finite support iteration in $K$ and hence Cohen reals are added at limit stages thus ${\rm cov}(\mathscr{M})=\lambda$ in $V[G]$.
By Theorem \ref{thmaddm} we have $\mathfrak{b}=\lambda$ in $V[G]$, so ${\rm add}(\mathscr{M})=\lambda$ due to Miller in \cite{MR613787}.
But Theorem \ref{thmaddm} also ensures $\clubsuit_{\aleph_1}$ in $V[G]$, so the proof is accomplished.

\hfill \qedref{coraddm}

There is an interesting upshot which follows from the above corollary.
We shall phrase and prove it in Corollary \ref{cortiltanfree} below, but we need some background.
Let $f$ be a function from $\kappa$ into $\mathcal{P}(\kappa)$. We say that $f$ is $C(\lambda,\mu)$ iff $|\bigcap\{f(x):x\in T\}|<\mu$ for every $T\in[\kappa]^\lambda$. We say that $f$ is $\kappa$-reasonable iff $f$ is both $C(\kappa,\omega)$ and $C(\omega,\kappa)$.
A subset $A$ of $\kappa$ is $f$-free iff $x\notin f(y)$ whenever $\{x,y\}\subseteq A$. Without the assumption that $f$ is $\kappa$-reasonable, one can find functions $f:\kappa\rightarrow\mathcal{P}(\kappa)$ for which there are no infinite free sets (and even no two-element free set, see \cite{MR3310341}). However, the mere assumption of $\kappa$-reasonability is insufficient for getting infinite free sets, as shown in the above mentioned paper.

It has been proved in \cite{MR3571281} that if $f$ is $\kappa$-reasonable and $\binom{\kappa}{\omega} \rightarrow \binom{\kappa}{\omega}_2$ then there exists an infinite free set for $f$. A natural question is whether this assumption is necessary. We shall give a negative answer to this question.
Recall that $\binom{\lambda}{\kappa} \rightarrow \binom{\lambda}{\kappa}_2$ holds iff for every coloring $c : \lambda \times \kappa \rightarrow 2$ there are $A \subseteq \lambda$ and $B \subseteq \kappa$ such that $|A|=\lambda, |B|=\kappa$ and $c \upharpoonright (A \times B)$ is constant.
Our first statement says that one can replace the assumption $\binom{\kappa}{\omega} \rightarrow \binom{\kappa}{\omega}_2$ by the weaker assumption $\binom{\kappa}{\kappa} \rightarrow \binom{\kappa}{\omega}_2$.

\begin{theorem}
\label{thmtiltanb} Infinite free subsets. \newline
Assume that $\binom{\kappa}{\kappa} \rightarrow \binom{\kappa}{\omega}_2$ and $f:\kappa\rightarrow \mathcal{P}(\kappa)$ is $\kappa$-reasonable. \newline
Then there exists an infinite free subset for $f$.
\end{theorem}

\par\noindent\emph{Proof}. \newline
We define a coloring $c: [\kappa]^2\rightarrow 2$ as follows. $c(\{\alpha,\beta\})=1$ iff $\alpha\notin f(\beta)\wedge \beta\notin f(\alpha)$. We employ the Erd\"os-Dushnik-Miller theorem to get either $H_0\in[\kappa]^\kappa$ such that $c\upharpoonright[H_0]^2=\{0\}$ or $H_1\in[\kappa]^\omega$ such that $c\upharpoonright[H_1]^2=\{1\}$. If there exists such $H_1$ then we are done, since it would be a free set for $f$ by the definition of the coloring $c$, so assume towards contradiction that there is no $H_1$ as above.

It follows that there is a $0$-monochromatic set $H_0$ of size $\kappa$.
We decompose it into $H_0=A\cup B$ such that $A\cap B=\emptyset$ and $|A|=|B|=\kappa$. Now we separate the cartesian product $A\times B$ into two disjoint collections:

$$
A\times B=\{\langle a,b\rangle: a\in f(b)\}\bigcup \{\langle a,b\rangle: a\notin f(b)\}.
$$

By the assumption $\binom{\kappa}{\kappa} \rightarrow \binom{\kappa}{\omega}^{1,1}_2$ we can choose $A_0\in[A]^\omega,B_0\in[B]^\kappa$ such that either $A_0\times B_0\subseteq\{\langle a,b\rangle: a\in f(b)\}$ or $A_0\times B_0\subseteq\{\langle a,b\rangle: a\notin f(b)\}$.

If $A_0\times B_0\subseteq\{\langle a,b\rangle: a\in f(b)\}$ then $A_0\subseteq \bigcap\{f(b):b\in B_0\}$, contradicting the assumption that $f$ is $C(\kappa,\omega)$. Similarly, if $A_0\times B_0\subseteq\{\langle a,b\rangle: a\notin f(b)\}$ then $b\in f(a)$ for every $a\in A_0,b\in B_0$ (since $a\notin f(b)$ and all the members are taken from $H_0$), so $B_0\subseteq \bigcap\{f(a):a\in A_0\}$, contradicting the assumption that $f$ is $C(\omega,\kappa)$.

\hfill \qedref{thmtiltanb}

Question 2.8 from \cite{MR3571281} is whether the assumption $\binom{\kappa}{\omega} \rightarrow \binom{\kappa}{\omega}_2$ is necessary for proving that every $\kappa$-reasonable function has an infinite free set.
In order to give a negative answer we will show that the positive relation $\binom{\kappa}{\kappa} \rightarrow \binom{\kappa}{\omega}_2$ is consistent with the negative relation $\binom{\kappa}{\omega} \nrightarrow \binom{\kappa}{\omega}_2$.

\begin{corollary}
\label{cortiltanfree} Let $\kappa=\cf(\kappa)>\aleph_1$. \newline
The following statements are consistent simultaneously:
\begin{enumerate}
\item [$(a)$] $2^\omega=\kappa$.
\item [$(b)$] $\binom{\kappa}{\omega} \nrightarrow \binom{\kappa}{\omega}_2$.
\item [$(c)$] tiltan holds.
\item [$(d)$] every $\kappa$-reasonable function $f:\kappa\rightarrow\mathcal{P}(\kappa)$ has an infinite free set.
\end{enumerate}
\end{corollary}

\par\noindent\emph{Proof}. \newline
Fix any $\kappa=\cf(\kappa)>\aleph_1$.
Using Theorem \ref{thmtiltanb} we force tiltan with $\mathfrak{b}=\mathfrak{c}=\kappa$.
It follows that $\mathfrak{b}=\mathfrak{d}=\kappa$ and hence $\binom{\kappa}{\omega} \nrightarrow \binom{\kappa}{\omega}_2$, see Theorem 2.2 of \cite{MR4068775}.
On the other hand, tiltan implies stick and by Proposition 3.3 of \cite{MR4068775} we have $\binom{\kappa}{\omega_1} \rightarrow \binom{\kappa}{\omega}_2$ and a fortiori $\binom{\kappa}{\kappa} \rightarrow \binom{\kappa}{\omega}_2$.
From Theorem \ref{thmtiltanb} we infer that every $\kappa$-reasonable function $f:\kappa\rightarrow\mathcal{P}(\kappa)$ has an infinite free set, so we are done.

\hfill \qedref{cortiltanfree}

We conclude this section with two problems.
The first one is about the relationship between superclub and ${\rm add}(\mathscr{M})$.
We do not know whether club can be replaced by superclub in the main result of this section.
The essential property that we used in the previous section in order to guarantee the preservation of superclub sequences seems problematic in the context of Hechler forcing.

\begin{question}
\label{qsuperadd} Is it consistent that superclub holds and ${\rm add}(\mathscr{M})>\omega_1$?
\end{question}

As mentioned at the beginning of this section, club and hence superclub imply ${\rm add}(\mathscr{N})=\aleph_1$, so the results of this section show that this cannot be stressed further in the horizontal dimension of Cicho\'ns diagram.
Brendle proved in \cite{MR2279653} that club is consistent with ${\rm cov}(\mathscr{N})=\aleph_2$, so in some sense we are covered in the vertical dimension as well.
In fact, it follows from the paper of Brendle that tiltan is consistent with an arbitrarily large value of ${\rm cov}(\mathscr{N})$, see \cite[Section 5]{MR2279653}.

\newpage

\section{Independent subsets and complete subgraphs}

Suppose that $\kappa,\mu\leq\lambda$.
The negative square brackets relation $\lambda\nrightarrow[\kappa]^2_\mu$ is the assertion that one can find a coloring $c:[\lambda]^2\rightarrow\mu$ such that for every $A\in[\lambda]^\kappa$ one has $c''[A]^2=\mu$.
Namely, every subset of $\lambda$ of size $\kappa$ assumes all the colors.
In the language of graph theory one can describe this relation as follows.
If $G=(V,E)$ is the complete graph with $\lambda$ vertices and one defines a decomposition $(E_\alpha:\alpha\in\mu)$ of $E$ then every subset $A\subseteq V$ of size $\kappa$ has at least one edge in every $E_\alpha$, that is $[A]^2\cap E_\alpha\neq\varnothing$ for every $\alpha\in\mu$.

The fact that $G$ is the complete graph is helpful here, since there are many edges and hence the probability that each $A$ of size $\kappa$ will intersect any $E_\alpha$ increases.
A graph theoretic generalization of this relation would be the same statement with respect to a larger class of graphs.
One has, however, to insert sufficiently many edges into $G$.
Recall that $I\subseteq V$ is \emph{independent} iff $I$ is edge-free, that is $[I]^2\cap E=\varnothing$.
A necessary requirement for our generalization is the absence of independent sets of size $\kappa$.

\begin{definition}
\label{defqproperty} $Q(\lambda,\kappa,\mu)$.
Assume that $\kappa,\mu\leq\lambda$. \newline
The property $Q(\lambda,\kappa,\mu)$ means that for every graph $G=(\lambda,E)$ with no independent subset of size $\kappa$ there exists a partition $(E_\alpha:\alpha\in\mu)$ of $E$ such that for every $\alpha\in\mu$ the graph $G_\alpha=(\lambda,E_\alpha)$ has no independent subset of size $\kappa$.
\end{definition}

The definition of $Q(\lambda,\kappa,\mu)$ comes from \cite{MR1387356}, and the fact that it is a strengthening of the square brackets relation with the same parameters is mentioned there.
Let us spell-out the proof.

\begin{obs}
\label{obsghk} If $Q(\lambda,\kappa,\mu)$ then $\lambda\nrightarrow[\kappa]^2_\mu$.
\end{obs}

\par\noindent\emph{Proof}. \newline
Let $G$ be the complete graph on $\lambda$ vertices and let $\{E_\alpha:\alpha\in\mu\}$ be a decomposition of $E$ ensured by the assumption $Q(\lambda,\kappa,\mu)$.
We define a coloring $c:[\lambda]^2\rightarrow\mu$ as follows:
$$
c(\gamma,\delta)=\alpha\Leftrightarrow\{\gamma,\delta\}\in E_\alpha.
$$
Assume that $A\in[\lambda]^\kappa$ and $\xi\in\mu$.
For every $\alpha\in\mu$ let $G_\alpha=(\lambda,E_\alpha)$.
By the assumption $Q(\lambda,\kappa,\mu)$ we see that $[A]^2\cap E_\alpha\neq\varnothing$ for every $\alpha\in\mu$, otherwise $A$ would be an independent subset of $G_\alpha$ of size $\kappa$.
Fix $\gamma,\delta\in A$ so that $\{\gamma,\delta\}\in E_\xi$.
By the definition of our coloring we see that $c(\gamma,\delta)=\xi$.
Since $\xi\in\mu$ was arbitrary we conclude that $c''[A]^2=\mu$ and the relation $\lambda\nrightarrow[\kappa]^2_\mu$ is therefore established.

\hfill \qedref{obsghk}

The property $Q(\lambda,\kappa,\mu)$ is strictly stronger than the relation $\lambda\nrightarrow[\kappa]^2_\mu$.
It has been proved in \cite{MR1387356} that the negation of $Q(\aleph_1,\aleph_1,\aleph_1)$ is consistent.
Contrariwise, Todor\v{c}evi\'c proved in \cite{MR908147} that $\aleph_1\nrightarrow[\aleph_1]^2_{\aleph_1}$ holds in \textsf{ZFC}.
Hence the implication in the above observation cannot be reversed.
However, for the negation of $Q(\aleph_1,\aleph_1,\aleph_1)$ one has to violate \textsf{CH}.
More generally, it is shown in \cite{MR1387356} that if $2^\kappa=\kappa^+$ then $Q(\kappa^+,\kappa^+,\kappa^+)$ holds.
Let us prove that this conclusion follows from a weaker assumption.
Namely, rather than the strong assumption $2^\kappa=\kappa^+$ one can assume the weaker assumption that superclub holds at $\kappa^+$ and still derive $Q(\kappa^+,\kappa^+,\kappa^+)$.

\begin{theorem}
\label{thmsuperclub} Let $\kappa$ be a regular cardinal.
If superclub holds at $\kappa^+$ then $Q(\kappa^+,\kappa^+,\kappa^+)$ holds.
\end{theorem}

\par\noindent\emph{Proof}. \newline
Suppose that $G=(\kappa^+,E)$ contains no independent subset of size $\kappa^+$.
Let $S$ be $S^{\kappa^+}_\kappa$, and let superclub be exemplified by $(S_\alpha:\alpha\in S)$.
Set $S_\alpha=\varnothing$ for every $\alpha\in\kappa^+-S$.
We would like to define a function $f:E\rightarrow\kappa^+$ with the following two properties at every $\alpha\in\kappa^+$:
\begin{enumerate}
\item [$(a)$] For every $\{\beta,\alpha\}\in E$ where $\beta<\alpha, f(\beta,\alpha)<\alpha$.
\item [$(b)$] If $\gamma<\alpha$ and $|\{\beta\in S_\gamma:\{\beta,\alpha\}\in E|=\kappa$ then $\{f(\beta,\alpha):\beta\in S_\gamma\wedge\{\beta,\alpha\}\in E\}=\alpha$.
\end{enumerate}
In order to define $f$ notice that if $\alpha\in\kappa^+$ then the cardinality of the set $x_\alpha=\{\gamma\in\alpha:|\{\beta\in S_\gamma:\{\beta,\alpha\}\in E|=\kappa\}$ is at most $|\alpha|\leq\kappa$.
So we choose a disjoint refinement $\{T_\gamma:\gamma\in x_\alpha\}$ for the collection $\{S_\gamma:\gamma\in x_\alpha\}$ and we fix a bijection $h^\alpha_\gamma:T_\gamma\rightarrow\alpha$.
For every $\gamma\in x_\alpha$ and every $\beta\in T_\gamma$ we let $f(\beta,\alpha)=h^\alpha_\gamma(\beta)$ and we define $f(\beta,\alpha)=0$ otherwise.

Having defined the function $f$ we can depict our decomposition of $E$.
For every $\tau\in\kappa^+$ set $\{\beta,\alpha\}\in E_\tau\Leftrightarrow f(\beta,\alpha)=\tau$.
Notice that $\{E_\tau:\tau\in\kappa^+\}$ is a decomposition of $E$ and we claim that it exemplifies $Q(\kappa^+,\kappa^+,\kappa^+)$ with respect to $G=(\kappa^+,E)$.

To see this, suppose that $H\subseteq\kappa^+$ and $|H|=\kappa^+$.
By induction on $\delta\in\kappa^+$ we choose an element $y_\delta\in H$ as follows.
As a first step, let $\{S_{\gamma_i}:i\in\kappa^+\}$ be a subset of the superclub sequence with the following two properties:
\begin{enumerate}
\item [$(a)$] $S_{\gamma_i}\subseteq H$ for every $i\in\kappa^+$.
\item [$(b)$] If $i<j<\kappa^+$ then $S_{\gamma_j}$ end-extends $S_{\gamma_i}$.
\end{enumerate}
Such a collection exists since one can choose $I\subseteq H, |I|=\kappa^+$ so that $S_I=\{\gamma\in\kappa^+:I\cap\gamma=S_\gamma\}$ is a stationary subset of $\kappa^+$.
Enumerate the elements of $S_I$ by $\{S_{\gamma_i}:i\in\kappa^+\}$ and verify the above properties.

Our second step is rendered by induction on $\delta\in\kappa^+$.
Suppose that $A_\delta=\{y_\zeta:\zeta\in\delta\}$ has been chosen and we try to choose $y_\delta$.
We will assume that $A_\delta\subseteq S_{\gamma_i}$ for some $i\in\kappa^+$, and let $i(\delta)$ be the first such ordinal.
Choose $y_\delta\in I\subseteq H$ so that $y_\delta>\gamma_{i(\delta)}$ and hence, in particular, $y_\delta>y_\zeta$ for every $\zeta\in\delta$ since $A_\delta\subseteq S_{\gamma_{i(\delta)}}\subseteq\gamma_{i(\delta)}$.
Notice that $A_{\delta+1}\subseteq S_{\gamma_j}$ for some $j\in\kappa^+$.

By way of contradiction assume that $\xi\in\kappa^+$ and $H$ is an independent subset of the graph $(\kappa^+,E_\xi)$.
This means that $\xi\notin f''[H]^2$.
By the construction, if $\delta>\xi$ then the set $z_\delta=\{\beta\in\delta:\{y_\beta,y_\delta\}\in E\}$ is of size less than $\kappa$ since otherwise the image of $f$ on these pairs will cover all the ordinals of $\delta$ including $\xi$.

We employ now the free subset theorem, which says (in this case) that if $h:\kappa^+\rightarrow\mathcal{P}(\kappa^+)$ satisfies $|h(\alpha)|<\kappa$ for every $\alpha\in\kappa^+$ then there exists an $h$-free subset $F\subseteq\kappa^+$ of size $\kappa^+$.
Consider the set-mapping $g:\kappa^+\rightarrow\mathcal{P}(\kappa^+)$ given by $g(\delta)=z_\delta$ for every $\delta\in\kappa^+$.
By the free subset theorem there is $F\subseteq\kappa^+,|F|=\kappa^+$ which is $g$-free.
Namely, if $\delta,\eta\in F$ and $\delta<\eta$ then $\delta\notin z_\eta$ and hence $\{y_\delta,y_\eta\}\notin E$.
It follows that $\{y_\beta:\beta\in F\}$ is an independent subset of $G=(\kappa^+,E)$, a contradiction.

\hfill \qedref{thmsuperclub}

We do not know whether $\clubsuit_{\kappa^+}$ is consistent with the failure of the principle $Q(\kappa^+,\kappa^+,\kappa^+)$.
It is clear, however, that the strategy of \cite{MR1489304} cannot yield such a result, since it cannot separate club from superclub.
On the other hand, the method of \cite{MR1623206} is promising in this context.
The idea is to force some property at $\aleph_2$ while keeping $\clubsuit_S$ at $S=S^{\omega_2}_\omega$ and then to collapse $\aleph_1$.
The collapse secures the club principle at $\omega_1$ but destroys superclub.
This method is limited to properties which can be forced at $\aleph_2$ and one also has to make sure that the pertinent property is preserved by the collapse.

\begin{question}
\label{qq} Is it consistent that $\clubsuit_{\kappa^+}$ holds but $Q(\kappa^+,\kappa^+,\kappa^+)$ fails for some infinite cardinal $\kappa$?
\end{question}

Let us indicate that the only property which separates tiltan from superclub, as far as we know, is Galvin's property.
Thus, a positive answer to the above question would be interesting.

\newpage

\bibliographystyle{alpha}
\bibliography{arlist}

\end{document}